\documentclass[12pt]{article}
\usepackage{graphicx}
\usepackage{amsfonts,amssymb,amsthm,eucal,amsmath}
\usepackage{amsbsy}
\usepackage{amscd}
\usepackage{url}
\usepackage[inner=30mm, outer=30mm, textheight=225mm]{geometry}

\newtheorem{thm}{Theorem}[section]

\newtheorem*{main_thm}{Theorem \ref{veering has angle struct}}
\newtheorem{lemma}[thm]{Lemma}
\newtheorem{prop}[thm]{Proposition}
\newtheorem{cor}[thm]{Corollary}
\theoremstyle{definition}
\newtheorem{defn}[thm]{Definition}
\newtheorem{que}[thm]{Question}

\theoremstyle{remark}
\newtheorem{rmk}[thm]{Remark}

\newcommand{\co}{\colon\thinspace}

\def\tri{\mathcal{T}}

\title{Veering triangulations admit strict angle structures}
\author{Craig D.\thinspace Hodgson, J.\thinspace Hyam Rubinstein, Henry Segerman\\ and Stephan Tillmann}

\begin{document}
\maketitle
\begin{abstract}
Agol recently introduced the concept of a veering taut triangulation, which is a taut triangulation with some extra combinatorial structure. We define the weaker notion of a ``veering triangulation'' and use it to show that all veering triangulations admit strict angle structures. We also answer a question of Agol, giving an example of a veering taut triangulation that is not layered.
\end{abstract}

\section{Introduction}

A basic question in 3-dimensional topology is to relate the combinatorics
of a triangulation of a 3-manifold to the geometry of the manifold.
The work of Gu\'eritaud and Futer~\cite{guer} deals with the case of Êhyperbolic structures
on once-punctured torus bundles and complements of two-bridge knots and links.
They study the angle structures on the natural layered ideal triangulations of these manifolds,
and use the volume maximization approach of Casson and Rivin Êto show that these
triangulations are geometric, i.e.\thinspace realized by positively oriented ideal hyperbolic tetrahedra.\\

In this paper, we introduce a new class of ``veering triangulations'' 
which includes the veering taut triangulations of Agol~\cite{agol_veering} and, in particular,  
the layered triangulations of once-punctured torus bundles.  Our main result 
shows that these veering triangulations admit strict angle structures.\\

Throughout this section, $M$ will denote the interior of an orientable 3--manifold with boundary a disjoint union of tori, 
imbued with a fixed ideal triangulation $\tri.$

\begin{defn}\label{taut_angle_structure}
An {\bf angle-taut tetrahedron} is an ideal tetrahedron equipped with an assignment of angles taken from $\{0,\pi\}$ to its edges so that two opposite edges are assigned $\pi$ and the other four are assigned $0$. A {\bf taut angle structure} on $M$ is an assignment of angles taken from $\{0,\pi\}$ to the edges of each tetrahedron in $M$ so that every tetrahedron is angle-taut and every edge of the triangulation has precisely two $\pi$ angles incident to it.
\end{defn}

\begin{figure}[htb]
\centering
\includegraphics[width=0.8\textwidth]{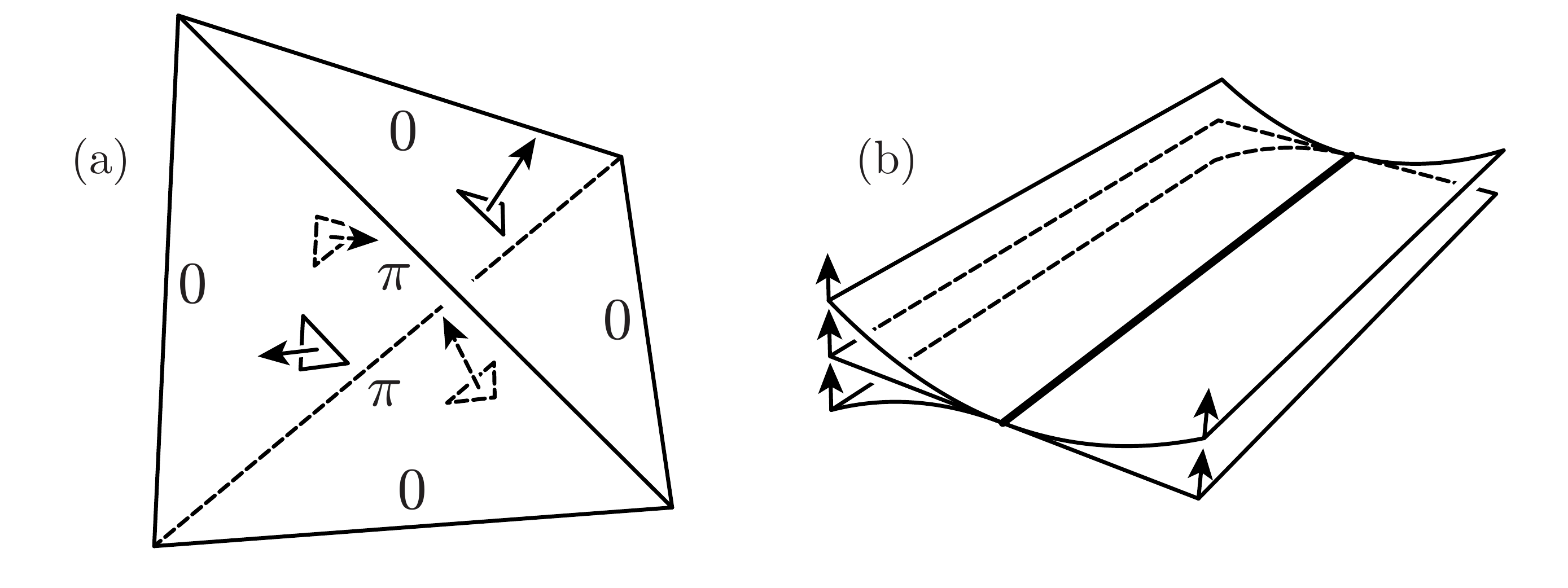}
\caption{Conditions for a taut ideal triangulation.}
\label{taut_ideal.pdf}
\end{figure}

The choice of adjective \emph{taut} for these angle structures is standard, but slightly unfortunate. Marc Lackenby \cite{lackenby00} introduced the following notion of a taut ideal triangulation in analogy with a taut foliation:

\begin{defn}\label{taut_structure}
A {\bf taut tetrahedron} is a tetrahedron such that each face is assigned a coorientation with two faces pointing inwards and two pointing outwards. Each edge of a taut tetrahedron is assigned an angle of either $\pi$ if the coorientations on the adjacent faces agree, or $0$ if they disagree. See Figure \ref{taut_ideal.pdf}(a) for the only possible configuration. Then $\tri$ is a {\bf taut ideal triangulation} of $M$ if there is a coorientation assigned to each ideal triangle so that every ideal tetrahedron is taut, and the sum of the angles around each edge is $2\pi$ (see Figure \ref{taut_ideal.pdf}(b)). This will also be called a {\bf taut structure} on $M.$
\end{defn}

A taut ideal triangulation comes with a compatible taut angle structure, but not every taut angle structure arises in this way. Examples of this are given in Section~\ref{non_layered_example}. \\

Ian Agol~\cite{agol_veering} recently introduced extra structure on a taut ideal triangulation, which he terms \emph{veering}, using a local condition. We define the notion of veering for any triangulation with a taut angle structure, and show that this new definition agrees with Agol's for taut ideal triangulations.

\begin{figure}[htb]
\centering
\includegraphics[width=0.25\textwidth]{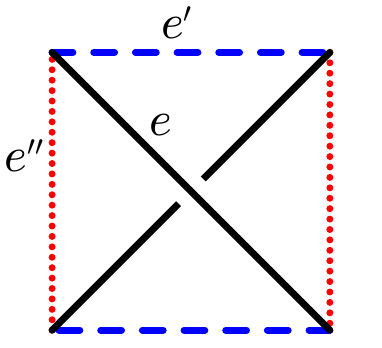}
\caption{The canonical picture of a veering tetrahedron. The $0$ angles are at the four sides of the square, and the $\pi$ angles are the diagonals. We indicate the veering directions on the $0$ angle edges of a tetrahedron by colouring the edges. Note that this picture depends on a choice of orientation for the tetrahedron.}
\label{veering_on_tetrahedron_d.pdf}
\end{figure}

\begin{defn}\label{veering_defn}
A {\bf veering tetrahedron} is an oriented angle-taut tetrahedron such that each edge with angle $0$ is coloured either red or blue (drawn dotted and dashed respectively)
as shown in Figure \ref{veering_on_tetrahedron_d.pdf}. We refer to the \underline{r}ed edges as {\bf right-veering} and the b\underline{l}ue edges as {\bf left-veering}. The $\pi$ angle edges are not assigned a colour with respect to the tetrahedron. A triangulation $\tri$ with a taut angle structure is a {\bf veering triangulation} of $M$ if there is a colour assigned to each edge in the triangulation so that every tetrahedron is veering. This will also be called a {\bf veering structure} on $M.$
\end{defn}

In this paper, angle-taut tetrahedra will generally be drawn flattened onto the page as in Figure \ref{veering_on_tetrahedron_d.pdf}, so that every tetrahedron appears as a quadrilateral with two diagonals, the four boundary edges all having angle $0$, and the two diagonals having angle $\pi.$ \\

In Section~\ref{sec:definitions}, we will interpret the definitions in terms of normal surface theory, and give a proof of the following result:

\begin{prop}\label{prop:equiv defns}
For a taut ideal triangulation, Definition \ref{veering_defn} is equivalent to the definition of a veering taut triangulation given by Agol \cite{agol_veering}.
\end{prop}

The main result of this paper is the following, a proof of which is given in Section~\ref{sec:veering implies angle}. Our more general notion of veering triangulations was motivated by this result, which does not depend on the existence of a taut triangulation of $M.$

\begin{thm}\label{veering has angle struct}
Veering triangulations of $3$-manifolds admit strict angle structures.
\end{thm}

In particular, using results of Casson and Thurston, it follows that any
$3$-manifold with a veering triangulation admits a complete hyperbolic structure.

\begin{que}Can all veering triangulations be realised as ideal hyperbolic triangulations in which all tetrahedra are positively oriented?
\end{que}

Veering triangulations seem to be very special (see Remark~\ref{rmk:veering data} for some data from the SnapPea census), but many manifolds have them. Agol \cite{agol_veering} proves that if we take any pseudo-Anosov mapping torus and puncture the surface along the singular points of the invariant measured foliations, then the mapping torus with the restricted monodromy has a {\bf layered triangulation} with compatible veering and taut structures (see \cite{agol_veering} for the sense of compatibility here). A layered triangulation is obtained by stacking tetrahedra on a triangulation of a surface. For instance, the canonical triangulations of once punctured torus bundles are all layered, taut and veering. Agol points out that the definition of veering does not depend on the triangulation being layered, and asks whether there is a veering taut triangulation which is not layered. In Section \ref{non_layered_example} we give such an example.

\begin{que}In addition to layering, what other ways are there to generate veering triangulations? 
\end{que}


\section{Definitions}
\label{sec:definitions}


\subsection{Ideal triangulation}

Let $M$ be a topologically finite 3-manifold, i.e.\thinspace the interior of a compact 3-manifold. An ideal triangulation $\tri$ of $M$ consists of a pairwise disjoint union of standard Euclidean 3--simplices, $\widetilde{\Delta} = \cup_{k=1}^{n} \widetilde{\Delta}_k,$ together with a collection $\Phi$ of Euclidean isometries between the 2--simplices in $\widetilde{\Delta},$ called {\bf face pairings}, such that $M \cong (\widetilde{\Delta} \setminus \widetilde{\Delta}^{(0)} )/ \Phi.$ The simplices in $M$ may be singular. It is well-known that every non-compact, topologically finite 3--manifold admits an ideal triangulation. We will assume throughout that $M$ as above is imbued with a fixed triangulation. We will also assume that $M$ is oriented and that all 
3--simplices in $M$ are coherently oriented.


\subsection{Quadrilateral types}

Let $\Delta^3$ be the standard 3--simplex with a chosen orientation. Suppose the edges from one vertex of $\Delta^3$ are labelled by $e,$ $e'$ and $e''$ so that the opposite edges have the same labelling. Then the cyclic order of $e,$ $e'$ and $e''$ viewed from each vertex depends only on the orientation of the 3--simplex, i.e.\thinspace is independent of the choice of vertex. It follows that, up to orientation preserving symmetries, there are two possible labellings, and we fix one of these labellings as shown in Figure~\ref{veering_on_tetrahedron_d.pdf}.\\

Each pair of opposite edges corresponds to a normal isotopy class of quadrilateral discs in $\Delta^3,$ disjoint from the pair of edges. We call such an isotopy class a {\bf normal quadrilateral type}. There is a natural cyclic order on the set of normal quadrilateral types induced by the cyclic order on the edges from a vertex, and this order is preserved by all orientation preserving symmetries of $\Delta^3.$ The particular cyclic order chosen corresponds to the 3--cycle $(q\ q' \ q''),$ where $q^{(k)}$ is dual to $e^{(k)}.$\\

Let $M^{(k)}$ be the set of $k$--simplices in $M$. If $\sigma \in M^{(3)},$ then there is an orientation preserving map $\Delta^3 \to \sigma$ taking the $k$--simplices in $\Delta^3$ to elements of $M^{(k)},$ and which is a bijection between the sets of normal quadrilateral types. This map induces a cyclic order of the normal quadrilateral types in $\sigma,$ and we denote the corresponding 3--cycle $\tau_\sigma.$ It follows from the above remarks that this order is independent of the choice of the map. \\

Let $\square$ denote the set of all normal quadrilateral types in $M.$ Define the permutation $\tau$ of $\square$ by
$$\tau = \prod_{\sigma \in M^3} \tau_\sigma.$$
See Figure~\ref{quad_types_in_veering_tet_d.pdf} for an illustration of the action of $\tau$ on the three quadrilateral types in a tetrahedron.\\

If $e\in M^{(1)}$ is any edge, then there is a sequence $(q_{n_1}, ..., q_{n_k})$ of normal quadrilateral types facing $e,$ which consists of all normal quadrilateral types dual to $e$ listed in sequence as one travels around $e.$ Then $k$ equals the degree of $e,$ and a normal quadrilateral type may appear at most twice in the sequence. This is well-defined up to cyclic permutations and reversing the order.


\subsection{Angle structures}

\begin{defn}\label{def:gen angle}

A function $\alpha\co \square \to \mathbb{R}$ is called a {\bf generalised angle structure on $M$} if it satisfies the following two properties:
\begin{enumerate}
\item If $\sigma^3 \in M^{(3)}$ and $q$ is a normal quadrilateral type supported by it, then
\begin{equation*}
   \alpha(q) + \alpha(\tau q) + \alpha(\tau^2 q) =\pi.
\end{equation*}
\item If $e\in M^{(1)}$ is any edge and $(q_{n_1}, ..., q_{n_k})$ is its normal quadrilateral type sequence, then
\begin{equation*}
\sum_{i=1}^k \alpha(q_{n_i})  =2\pi.
\end{equation*}
\end{enumerate}
\end{defn}

Dually, one can regard $\alpha$ as assigning angles $\alpha(q)$ to the two edges opposite $q$ in the tetrahedron containing $q$.\\

A generalised angle structure is called 
a {\bf taut angle structure on $M$} if its image is contained in $\{0,\pi\},$ 
a {\bf semi-angle structure on $M$} if its image is contained in $[0,\pi],$ and 
a {\bf strict angle structure on $M$} if its image is contained in $(0,\pi).$


\subsection{Agol's definition}

Throughout this subsection, suppose $M$ has a veering triangulation with underlying taut angle structure $\alpha \co \square \to \{0, \pi\}.$ The conventions regarding orientations and normal quadrilateral types immediately imply the following fact (compare Figure \ref{quad_types_in_veering_tet_d.pdf}):

\begin{lemma}\label{lem:veering quad}
Suppose $\sigma$ is a tetrahedron in $M.$ Let $q$ be the quadrilateral type dual to the edges with label $\pi.$ Then the edges dual to $\tau(q)$ are left-veering (blue) and the edges dual to $\tau^2(q)$ are right-veering (red).
\end{lemma}

If $e\in M^{(1)}$ is any edge and $(q_{n_1}, ..., q_{n_k})$ is its normal quadrilateral type sequence, then there are precisely two normal quadrilateral types in this sequence on which $\alpha$ takes the value $\pi.$ These separate the sequence into two subsequences of consecutive 0--angle quadrilateral types, which we call the \textbf{sides} of $e.$ Their lengths are called the {\bf one-sided degrees} of $e.$

\begin{figure}[htb]
\centering
\includegraphics[width=0.35\textwidth]{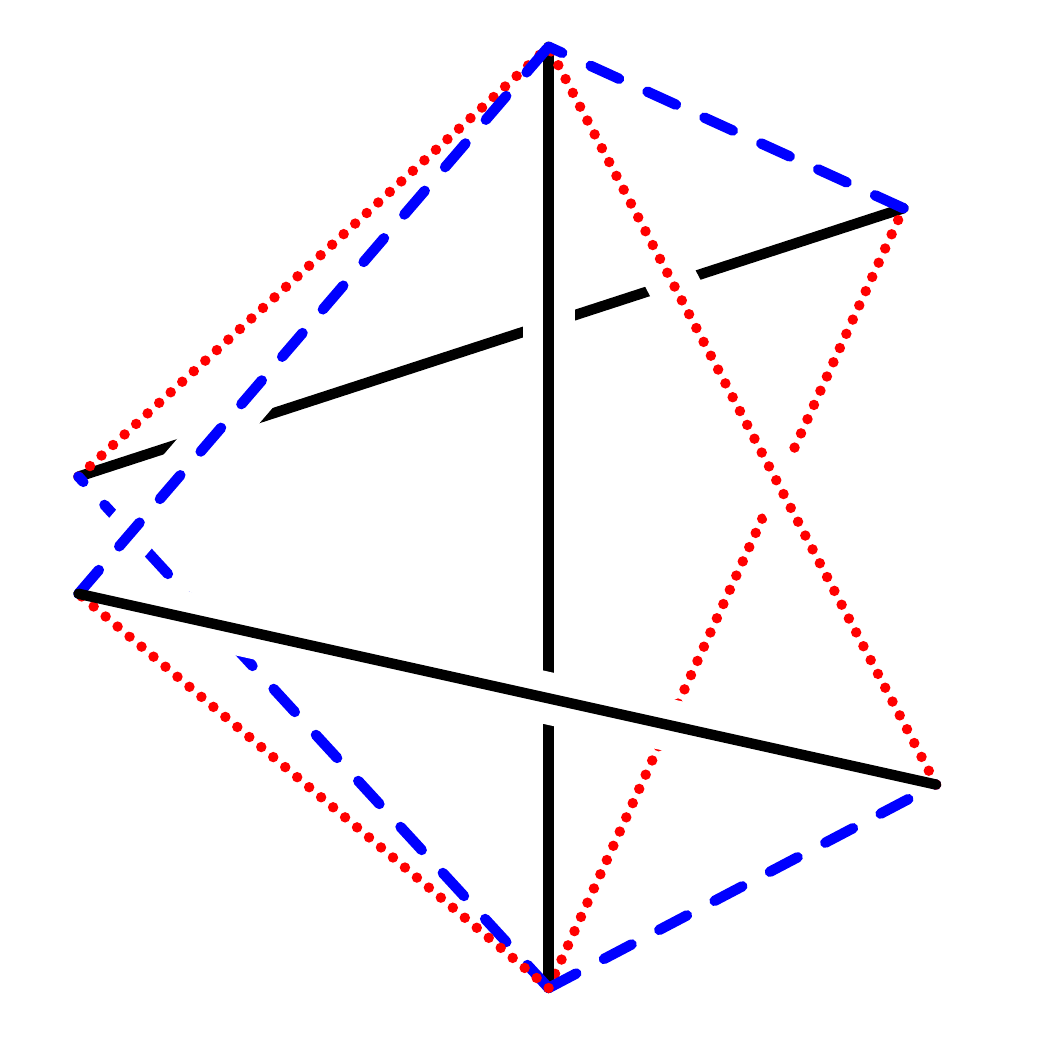}
\caption{The one-sided degree in a veering taut triangulation is at least one.}
\label{one_sided_degree_at_least_1_d.pdf}
\end{figure}

\begin{lemma}\label{lem:one-sided degree}
Each one-sided degree of $e$ is at least one.
\end{lemma}

\begin{proof}
At each vertex of a veering tetrahedron there is a cyclic order of the three edges of the tetrahedron, taking the $\pi$--edge to the left-veering (blue) edge to the right-veering (red) edge. Notice that the order is the same at every vertex. This implies that the two normal quadrilateral types dual to the edge $e$ having angle $\pi$ cannot be adjacent in the normal quadrilateral type sequence of $e,$ since otherwise there would be a conflict in the colouring. See Figure~\ref{one_sided_degree_at_least_1_d.pdf}.
\end{proof}

\begin{cor}
The degree of each edge is at least four.
\end{cor}

We now show that for a taut ideal triangulation, Definition \ref{veering_defn} is equivalent to the definition of a veering taut triangulation given by Agol (\cite{agol_veering}, Definition 4.1). According to that definition, one needs to check two conditions for each edge in the triangulation. The second condition is already given by Lemma~\ref{lem:one-sided degree}. The first is verified by the following lemma. See \cite{agol_veering} for further context of the terms involved in its statement. 

\begin{lemma}\label{veering_defn_original}
Suppose that the triangulation of $M$ is taut and veering, and has a taut angle structure compatible with both.
Let $e$ be a right (resp.\thinspace left) -veering edge in $M$. Consider the sequence of oriented triangles incident to one side of $e$, and the vertices of these triangles opposite $e$. If we order the vertices moving ``upwards'' according to the coorientation of the triangles, they are also ordered from left to right (resp.\thinspace right to left) as viewed from $e$. Similarly for the other side of $e$. 
\end{lemma}

\begin{figure}[htb]
\centering
\includegraphics[width=0.7\textwidth]{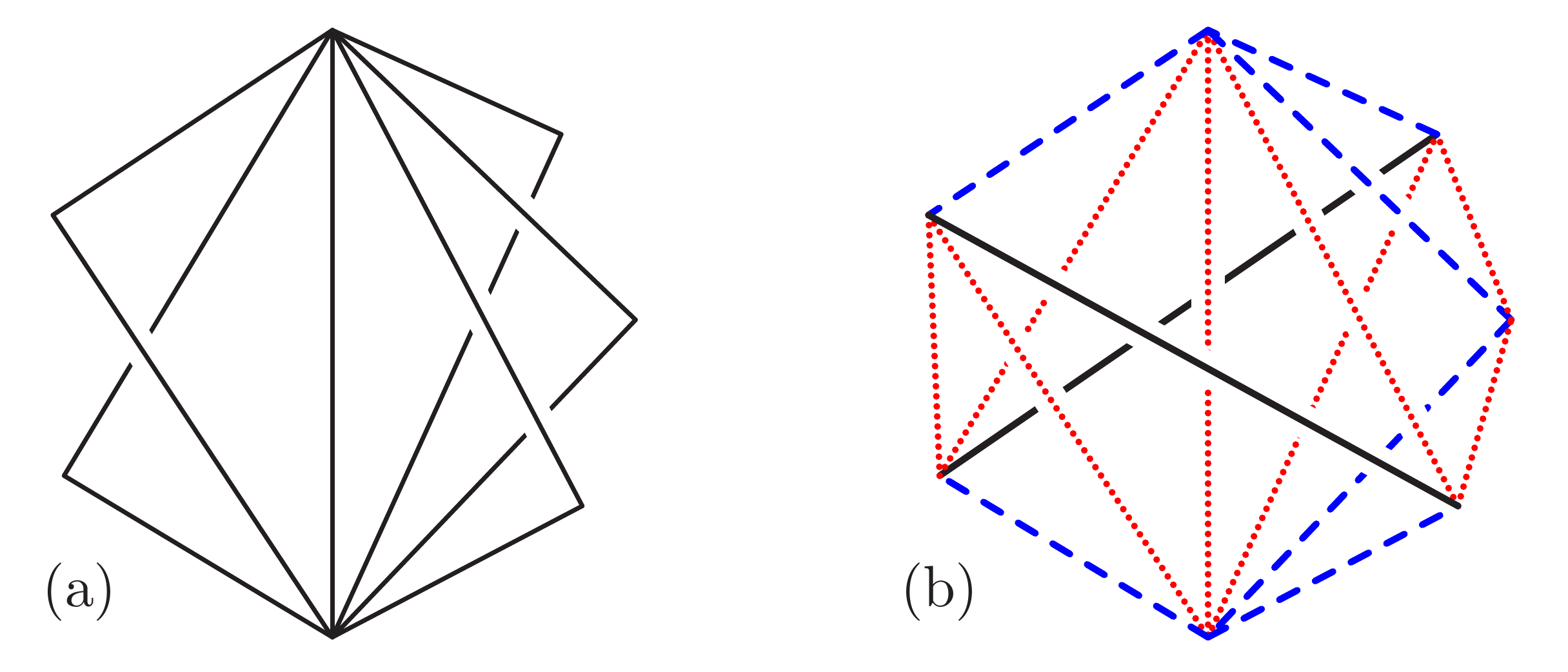}
\caption{The triangles adjacent to a right-veering edge, and veering directions on the link of a right-veering edge.}
\label{right_veering_triangles_and_link_d.pdf}
\end{figure}

See Figure \ref{right_veering_triangles_and_link_d.pdf}(a) for a picture of triangles incident to a right veering edge, veering to the right as the triangles are ordered approaching the reader. 

\begin{proof}
The result is a consequence of stacking copies of tetrahedra coloured as in Definition \ref{veering_defn} at $e$ so that the colours match. We have already observed 
the pattern on the two tetrahedra with label $\pi$ at $e$ in Figure~\ref{one_sided_degree_at_least_1_d.pdf}. The result now follows by induction on the number of tetrahedra abutting $e$ by considering the cases of whether $e$ is left-veering (blue) or right-veering (red). See Figure~\ref{right_veering_triangles_and_link_d.pdf}(b) for an example.
\end{proof}


\begin{proof}[Proof of Proposition~\ref{prop:equiv defns}]
Lemmas~\ref{lem:one-sided degree} and \ref{veering_defn_original} immediately imply the forward direction of Proposition~\ref{prop:equiv defns}. For the reverse direction: a taut triangulation satisfying the conclusion of Lemma \ref{veering_defn_original} has veering directions on the $0$ angle edges of each tetrahedron matching the colouring given by Definition \ref{veering_defn}. To see this, observe that each of the $0$ angle edges in Figure \ref{veering_on_tetrahedron_d.pdf} has a pair of incident triangles, determining the veering direction.
\end{proof}

\begin{rmk}
Lemma \ref{lem:veering quad} points to a further generalisation of veering triangulations, separated from an underlying taut angle structure. It consists merely of a choice of normal quadrilateral type $q$ in each oriented tetrahedron subject to compatibility of edge veering directions, where for each tetrahedron the left and right-veering edges are defined to be those dual to $\tau(q)$ and $\tau^2(q)$ respectively. In the current sense of veering, for each tetrahedron the choice of $q$ satisfies $\alpha(q)=\pi$.
\end{rmk}


\section{Veering triangulations and angle structures}
\label{sec:veering implies angle}

The main result of this paper is the following:

\begin{main_thm}
Veering triangulations of $3$-manifolds admit strict angle structures.
\end{main_thm}

Work of Kang and Rubinstein~\cite{kang_rubinstein_2005} and Luo and Tillmann~\cite{luo_tillmann_2008} links the existence of angle structures to the normal surface theory of $M$ using duality principles from linear programming. The {\bf normal surface solution space}  $C(M; \tri)$ is a vector subspace of $\mathbb{R}^{7n}$ where $n$ is the number of tetrahedra in $\tri$, consisting of vectors  satisfying the {\bf compatibility equations} of normal surface theory. The coordinates of $x \in \mathbb{R}^{7n}$ represent weights of the four normal triangle types and the three normal quadrilateral types in each tetrahedron, and the compatibility equations state that normal triangles and quadrilaterals have to meet the 2--simplices of $\tri$ with compatible weights. \\


A vector in $\mathbb{R}^{7n}$ 
is called {\bf admissible} if at most one quadrilateral coordinate from each tetrahedron is non-zero and all coordinates are non-negative. An integral admissible element of $C(M; \tri)$ corresponds to a unique embedded, closed normal surface in $M$ and vice versa. As a reference for other facts from normal surface theory, please consult~\cite{jaco_oertel}.\\

There is a linear function $\chi^*\co C(M; \tri) \to \mathbb{R}$ which agrees with the Euler characteristic $\chi$ on embedded and immersed normal surfaces. 
If $\tri$ admits a generalised angle structure $\alpha$, then the formal Euler characteristic $\chi^*$ can be computed by
%
$$2\pi \chi^*(x)= \sum_q -2\alpha(q) x_q,$$
where $x_q$ is the normal coordinate of the normal quadrilateral type $q$.\\

Since a taut angle structure has image $\{0, \pi\},$ we have the following simple  application of Theorems 1 and 3 of \cite{luo_tillmann_2008}:

\begin{cor}\label{luo_tillmann}
Suppose $M$ has a taut angle structure. Then the following are equivalent
\begin{enumerate}
\item $M$ admits a strict angle structure.
\item For all $x \in C(M; \tri)$ with all quadrilateral coordinates non-negative and at least one quadrilateral coordinate positive, $\chi^*(x) < 0.$
\item There is no $x \in C(M; \tri)$ with all quadrilateral coordinates non-negative and at least one quadrilateral coordinate positive and $\chi^*(x) = 0.$
\end{enumerate}
\end{cor}

\begin{figure}[htb]
\centering
\includegraphics[width=0.82\textwidth]{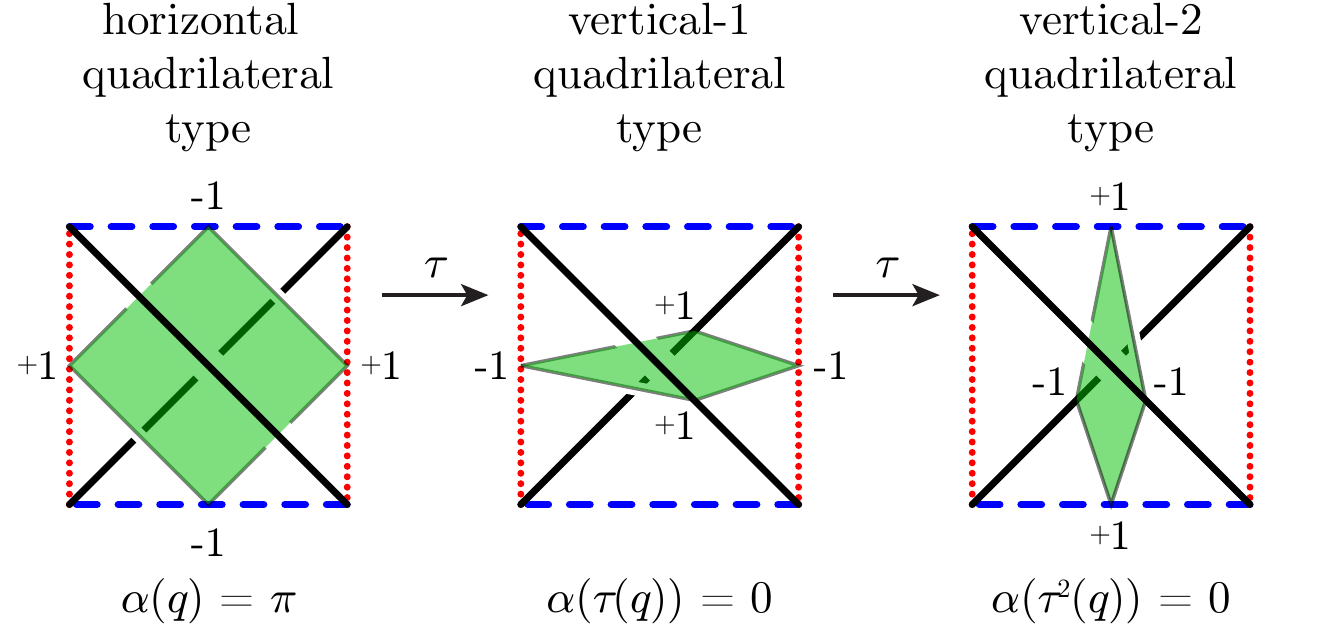}
\caption{The one horizontal and two vertical quadrilateral types in an angle-taut tetrahedron, with associated signs in the $Q$-matching equations.}
\label{quad_types_in_veering_tet_d.pdf}
\end{figure}

The quadrilateral coordinates in the solutions to the compatibility equations satisfy the so-called {\bf Q--matching equations} due to Tollefson \cite{Tollefson}. Given the sequence 
$(q_{n_1}, ..., q_{n_k})$ of normal quadrilateral types facing an edge $e,$ one associates a sign $\varepsilon(q)\in \{\pm 1\}$ to each element $q$ in the sequence  $(\tau(q_{n_1}), \tau^2(q_{n_1}), ..., \tau(q_{n_k}), \tau^2(q_{n_k}))$ of quadrilateral types incident with $e.$ If a normal quadrilateral type appears more than once in this sequence, it may have different signs. The $Q$--matching equation associated to $e$ is
$$\sum \varepsilon(q) x_q = 0,$$
where the sum is taken over the elements of $(\tau(q_{n_1}), \tau^2(q_{n_1}), ...,  \tau(q_{n_k}), \tau^2(q_{n_k})).$ In the situation of a veering triangulation, the signs are as given in Figure~\ref{quad_types_in_veering_tet_d.pdf}. This suggests another alternate interpretation of the colouring on a veering tetrahedron in terms of the signs of the quadrilateral type dual to the $\pi$ edges.

\begin{defn}
For a triangulation with a taut angle structure the three quadrilateral types  within each tetrahedron fall into two classes:
(i) there are two kinds of {\bf vertical} quadrilateral types (denoted {\bf vertical-1} and {\bf vertical-2}) which have two $0$ angle and two $\pi$ angle corners,  and
(ii) one kind of {\bf horizontal} quadrilateral type, which has four $0$ angle corners.
If $q$ is the horizontal quadrilateral type in a tetrahedron, then $\tau(q)$ and $\tau^2(q)$ are the vertical-1 and vertical-2 quadrilateral types respectively. See Figure~\ref{quad_types_in_veering_tet_d.pdf}.
\end{defn}

Note that a positive coordinate corresponding to a horizontal quadrilateral type contributes negatively to $\chi^*(x)$, whereas those corresponding to vertical quadrilateral types contribute $0$.  We remark that the set of solutions to the $Q$--matching equations also includes the so-called spun-normal surfaces, which do not appear as solutions in $C(M; \tri).$  Thus non-existence of solutions to the $Q$--matching equations under the conditions of Corollary~\ref{luo_tillmann}(3) implies (but is not implied by) non-existence of $x \in C(M; \tri)$ with those conditions, and we have the following:

\begin{lemma}\label{no_vert_solns}
A triangulation with a taut angle structure admits a strict angle structure if it admits no non-negative solution to the $Q$-matching equations with at least one quadrilateral coordinate positive and for which every non-zero quadrilateral type is vertical. 
\end{lemma} 

Hence Theorem \ref{veering has angle struct} follows from the following:

\begin{prop}
If $M$ has a veering triangulation, then there is no non-negative solution to the $Q$-matching equations with at least one quadrilateral coordinate positive and for which every non-zero quadrilateral is vertical.
\end{prop}

\begin{proof}
By way of contradiction, suppose there is a non-negative solution $x$ to the $Q$-matching equations with at least one quadrilateral coordinate positive and for which every non-zero quadrilateral type is vertical. \\

Suppose $q$ is a vertical-1 quadrilateral type with $x_q>0,$ supported by the tetrahedron $\sigma$. Then the signs associated to its corners on the red edges of $\sigma$ are both negative, and hence its contribution to the respective $Q$--matching equations is negative. This cannot be compensated for by vertical-2 quadrilateral types, since their positive signs are on blue edges. It follows that it must be compensated for by vertical-1 quadrilateral types.\\

This implies that the total negative contribution of all non-zero vertical-1 quadrilateral types to the sum of all $Q$--matching equations of the red edges equals the total positive contribution of all vertical-1 quadrilateral types to the sum of all $Q$--matching equations of the red edges. Hence a vertical-1 quadrilateral type with non-zero coordinate must have all four corners on red edges; similarly for vertical-2 quadrilateral types and blue edges. So the $Q$--matching equations restricted to red edges and vertical-1 quadrilateral types are independent of the $Q$--matching equations restricted to blue edges and vertical-2 quadrilateral types. Therefore we may assume that $x_q=0$ for each vertical-2 quadrilateral type. This implies that the solution $x$ is admissible. Since $x$ is a convex linear combination of the so-called admissible vertex solutions, we may assume that $x$ is an admissible vertex solution, and may hence choose a connected, oriented normal surface $S$ in the projective normal isotopy class defined by $x.$ Since $S$ is oriented, we can choose a transverse orientation to $S.$\\

The argument in the previous two paragraphs goes through verbatim if we interchange ``vertical-1'' with ``vertical-2'' and ``red'' with ``blue'', as does the rest of the proof (also interchanging ``$+1$'' with ``$-1$'').\\

We now analyse how the quadrilateral discs of the surface $S$ sit in the triangulation. First notice that each quadrilateral disc has all of its corners on red edges. The two corners with sign $+1$ inherit angle $\pi,$ and the corners with sign $-1$ inherit angle $0.$ Since $S$ is embedded and the sum of angles around any edge equals $2\pi,$ it follows that at any vertex in the induced cell structure of $S,$ there are at most four quadrilateral discs.\\

The next claim pertains to the tetrahedra with four edges of the same colour. As noted above, all tetrahedra containing quadrilaterals of $S$ have this property.\\ 

\emph{Property ($*$): Suppose edge $e$ is red, and has the sequence $(q_{n_1}, \ldots q_{n_k})$ of quadrilateral types  facing it, where $\alpha(q_{n_1})=0$ and the tetrahedron supporting $q_{n_1}$ has four red edges. Then $(q_{n_1})$ is a side of $e.$ In other words, $\alpha(q_{n_2})=\alpha(q_{n_k})=\pi.$}\\

We now prove Property ($*$). Consider a tetrahedron $\sigma$ with both $\pi$ angles right-veering. See Figure \ref{tet_with_4_red_d.pdf}. Then $\sigma$ has angle $0$ at the edge $e$. The veering choices on $T$ force that the two (possibly non-distinct) tetrahedra adjacent to the two faces of $T$ incident to $e$ have $\pi$ angles at $e$. This is illustrated in the right diagram of Figure \ref{tet_with_4_red_d.pdf}.

\begin{figure}[htb]
\centering
\includegraphics[width=0.6\textwidth]{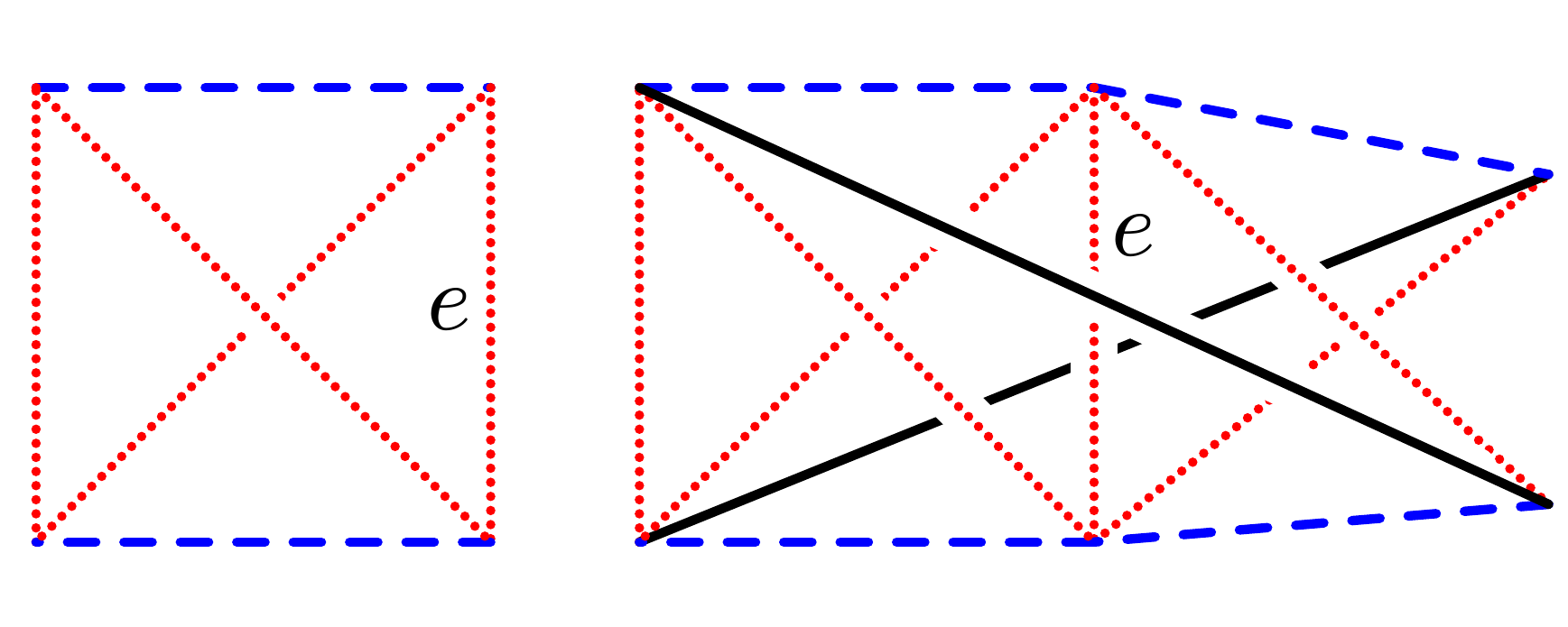}
\caption{On the left: a tetrahedron with both $\pi$ angles right-veering (red). On the right: this choice of veering implies that this tetrahedron contributes the only $0$ angle on its side of $e$. }
\label{tet_with_4_red_d.pdf}
\end{figure}

Property ($*$) places severe restrictions on the cell structure of $S.$ If a quadrilateral disc meets an edge $e$ with one of its $-1$ corner, then it is contained in a tetrahedron corresponding to the $(q_{n_1})$ side of $e.$ Thus, at each vertex in the induced cell structure of $S,$ either precisely two quadrilaterals or precisely four quadrilaterals meet. Moreover, the former is the case if $\deg(e)\ge 5,$ and the latter can only happen if $\deg(e)=4.$ In addition, each quadrilateral in $S$ meets quadrilaterals (rather than triangles) along at least two of its edges.\\

We will now use the transverse orientation and the natural geometry of quadrilaterals to start exploring the subsurface of $S$ made up of quadrilateral discs. See Figure \ref{traverse_quads.pdf}. Suppose we start on $q_0.$ We set our compass so that one of the positive corners is north, then the other is south and the negative corners are east and west. As we walk from $q_0$ to another quadrilateral disc, we parallel transport our compass, and do not keep track of whether we visited this disc previously. Given an arbitrary quadrilateral disc $q$, Property ($*$) implies that we can always find a quadrilateral disc in $S$ meeting $q$ along at least one of its two sides incident to a $-1$ corner. For $q_0$ this means that it meets a quadrilateral disc along either the north-east or south-east edge. This is true for any quadrilateral disc we have parallel transported our compass to.\\

\begin{figure}[htb]
\centering
\includegraphics[width=0.7\textwidth]{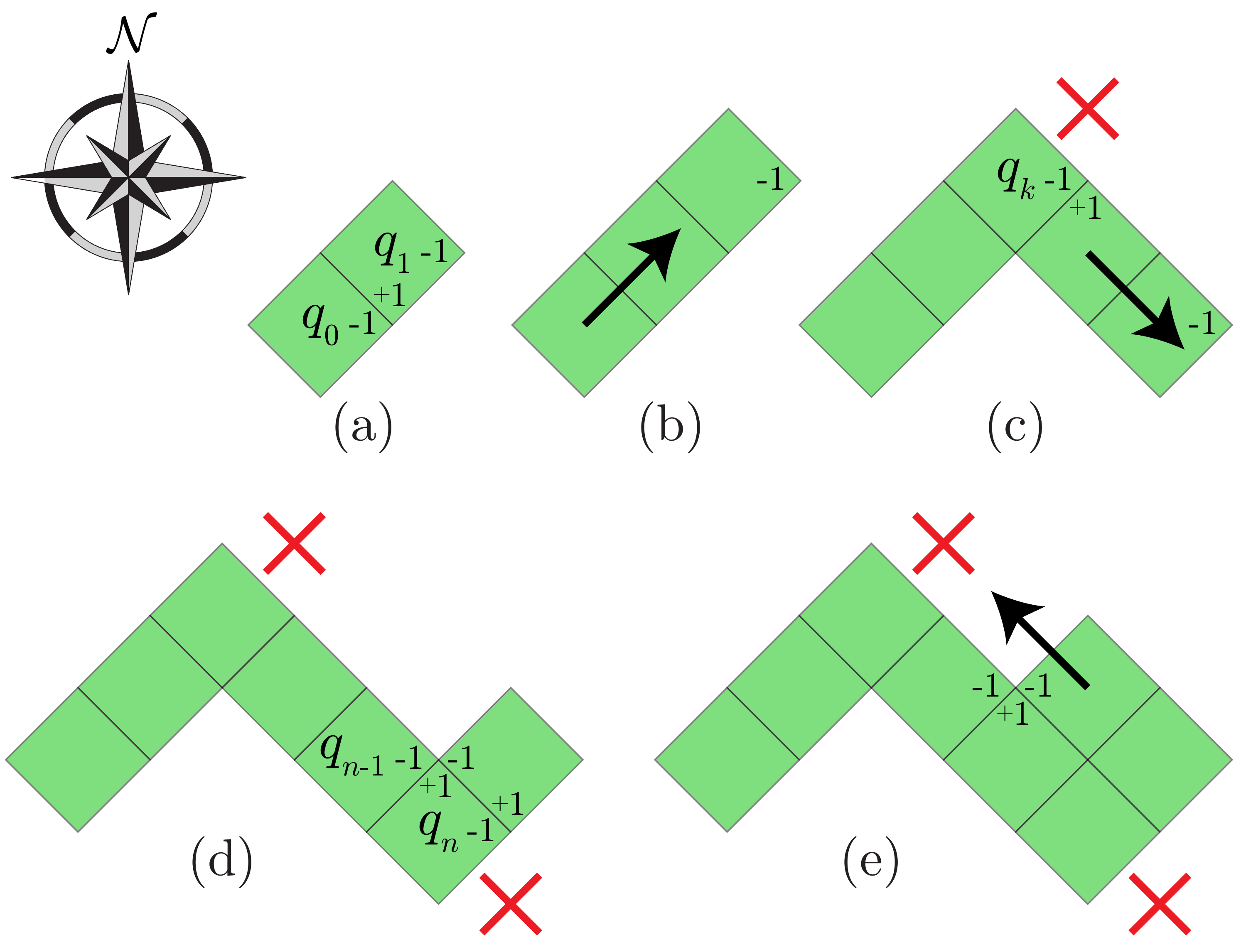}
\caption{Traversing through quadrilateral discs.}
\label{traverse_quads.pdf}
\end{figure}

Suppose we cannot walk from $q_0$ indefinitely in the north-easterly direction. Then after travelling in the north-easterly direction across finitely many quadrilateral discs, we reach a quadrilateral disc $q_k,$ which meets a normal triangle along its north-east side. Since the eastern corner has a $-1,$ $q_k$ must meet a quadrilateral disc along its south-east side.\\

\emph{Claim: We can walk from $q_k$ indefinitely in the south-easterly direction.}\\

The proof is by contradiction and sketched in Figure~\ref{traverse_quads.pdf}. Suppose to the contrary that after traversing discs $q_{k+1}, \ldots, q_{n-1},$ we reach a disc $q_n$ which does not meet another quadrilateral disc along its south-east side. Then it must meet another along its north-east side in order to cancel its eastern $-1$ corner. But then the edge at its northern $+1$ corner is incident with three quadrilaterals, and hence must be incident with four. This implies that this edge is of degree four. But then $q_{n-1}$ is glued along its north-east edge to another quadrilateral. This now propagates all the way to $q_k,$ giving a contradiction. This proves the claim.

\begin{figure}[htb]
\centering
\includegraphics[width=0.37\textwidth]{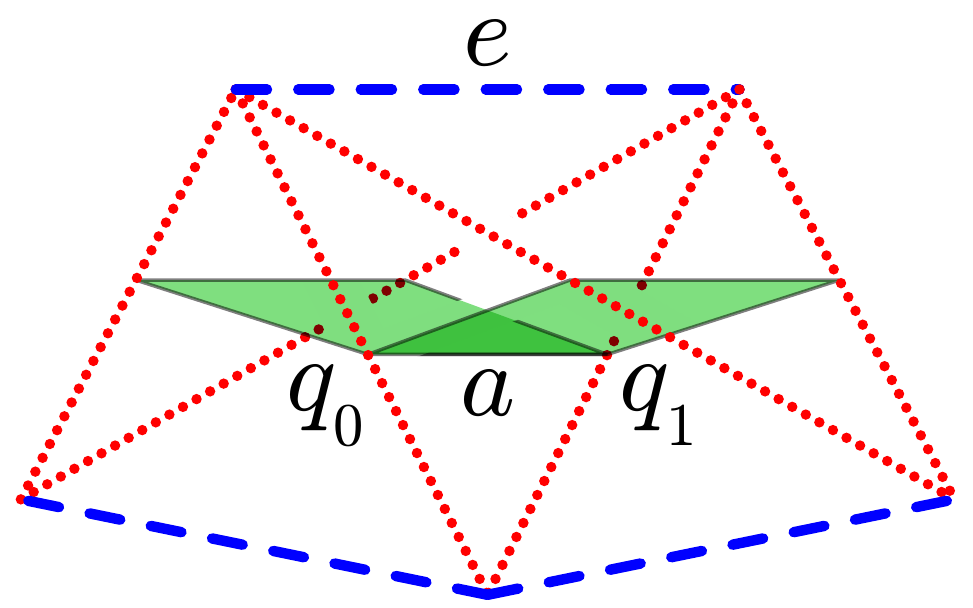}
\caption{The dual edge remains the same as we walk across quadrilateral discs.}
\label{tet_with_4_red_quads_d.pdf}
\end{figure}

In all cases then, there is a quadrilateral disc $q_0$ in $S,$ from which one can walk indefinitely across quadrilaterals. Let $a$ be the normal arc in the boundary of $q_0$ that we walk across (as shown in Figure~\ref{tet_with_4_red_quads_d.pdf}). We may assume that the transverse orientation of $S$ induced on $a$ points not to the ideal vertex of the triangle containing $a,$ but to the ideal edge $e$ which $a$ does not meet. Then $q_1$ again has its transverse orientation pointing towards $e$ and is also dual to it. It follows inductively that all quadrilateral discs in the indefinite path are dual to $e,$ and that all quadrilateral types dual to $e$ appear in $S$ and are therefore all vertical. But then the total angle sum at $e$ is zero, contradicting the fact that $M$ has a taut angle structure.
\end{proof}

\begin{rmk}
If the given manifold with veering triangulation is known to be hyperbolic, then we can alter the proof in the following way: Since the surface $S$ is embedded, we have $\chi(S) = \chi^*(x)=0,$ and since it is orientable and connected it must be a torus or a spun-normal annulus. Theorem 2.6 of \cite{kang_rubinstein_2005} then tells us that this torus must be non-boundary parallel and essential, leading to a contradiction to the assumption that the manifold is hyperbolic. If there is an annulus, there are two cases. If the annulus is boundary parallel, one can use this surface as a barrier (as in \cite{0-efficient}) and find a normal torus which is topologically but not normally parallel to the boundary, giving a contradiction as in \cite{kang_rubinstein_2005}. If there is an essential annulus then the manifold is again not hyperbolic. 
\end{rmk}


\section{A veering taut triangulation that is not layered}\label{non_layered_example}

Agol's construction~\cite{agol_veering} gives us the existence of many veering taut triangulations, but they are all layered triangulations. In this section we give an example of a veering taut triangulation that is not layered. In fact, the manifold does not fibre over the circle. We also list some examples of veering triangulations that do not fibre over the circle, but are not taut triangulations.\\

\begin{figure}[htb]
\centering
\includegraphics[width=0.9\textwidth]{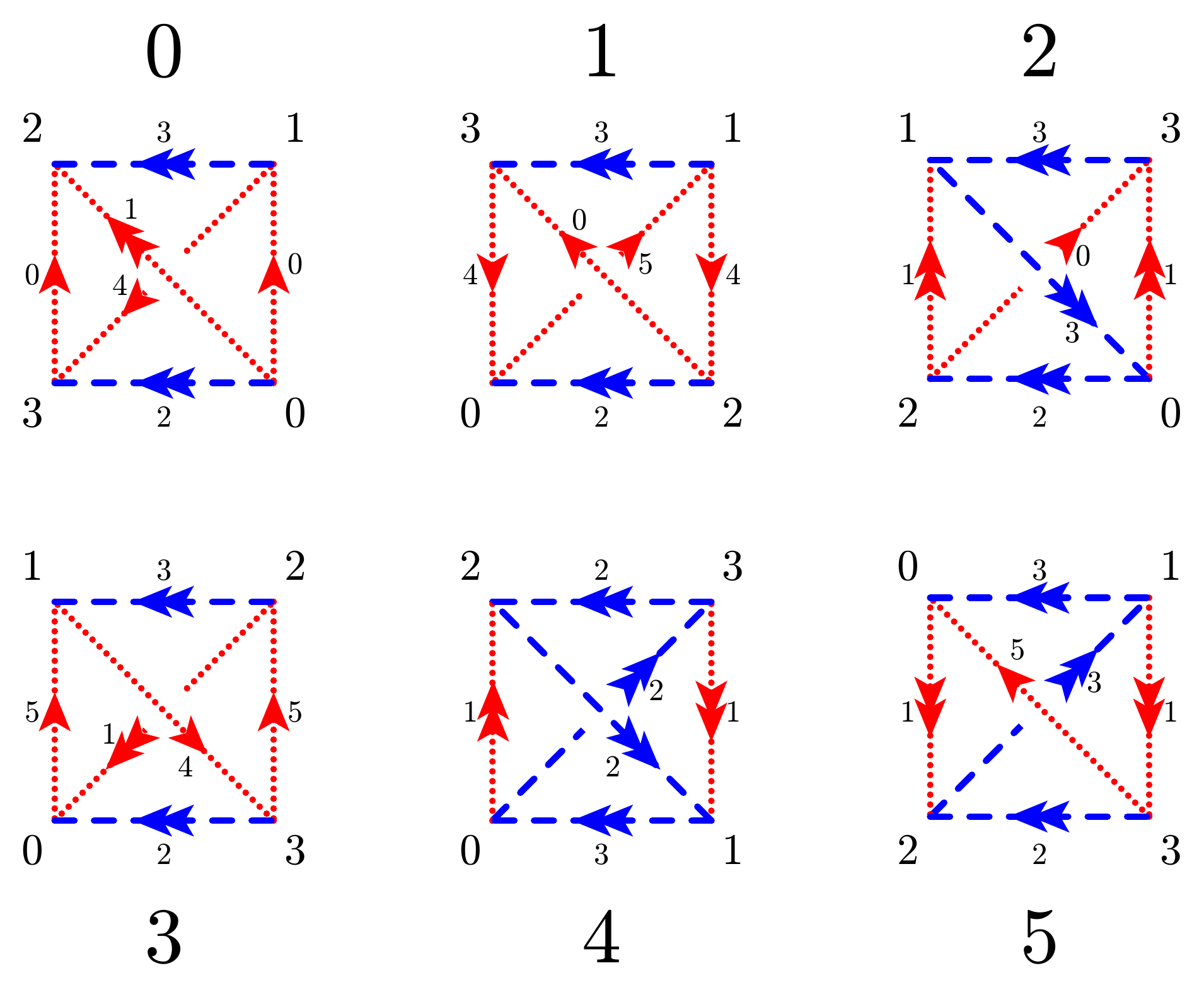}
\caption{The triangulation s227 from the SnapPea census, with the veering taut structure indicated. There are six tetrahedra labelled from 0 to 5 and six edges, also labelled from 0 to 5. Single arrows refer to degree 4 edges and double arrows degree 8 edges. We have drawn the tetrahedra so that they are taut, flattened onto the page as in Figure \ref{veering_on_tetrahedron_d.pdf}, with consistent coorientation of the triangles in the direction out of the page towards the reader. Moreover, the orientation of each tetrahedron agrees with the vertex sequence $0,1,2,3.$}
\label{s227_d.pdf}
\end{figure}

The main example is the manifold s227 from the SnapPea census~\cite{snappea}, with the triangulation as given in the census. See Figure \ref{s227_d.pdf}. One can show that the manifold s227 is not fibred by applying a method due to Brown \cite{brown1987} to compute the BNS invariant. More specifically, the fundamental group of the manifold is isomorphic to
$$\Gamma = \langle a, b \mid a^4b^2a^{-1}b^{-1}a^{-1}b^2a^{-1}b^{-1}a^{-1}b^2=1 \rangle,$$
and hence has two generators and one relator, with abelianisation isomorphic to $\mathbb{Z} \oplus \mathbb{Z}_4.$ An application of the algorithm described in Theorem 4.2 of \cite{brown1987} verifies that there is no epimorphism $\Gamma \to \mathbb{Z}$ with finitely generated kernel.

\begin{rmk}
This example was found by searching the SnapPea census using the Regina software package by Ben Burton~\cite{regina}, together with a listing of SnapPea census manifolds with data on whether or not they fibre compiled by Nathan Dunfield (\url{http://dunfield.info/snappea/tables/mflds_which_fiber}). There is no guarantee that this is the ``smallest'' example of a non-fibred veering taut triangulation, since we only checked the triangulation of each manifold as it appears in the census. In addition we have not checked for layering directly; there are a number of triangulations in the census with taut angle structures that are both taut and veering, and it is likely that some of these examples are not layered even if the manifold fibres.
\end{rmk}

\begin{rmk}
Examples of veering triangulations of manifolds that do not fibre over the circle but that are not taut triangulations were found by a similar method. They are s438, s772, s773, s779, v3128, v3243, v3244, v3377 and v3526, again with the triangulations as given in the SnapPea census and taut angle structure as given in Table \ref{table_example_taut_angles}. These manifolds do not fibre either by the BNS invariant again (as calculated by Dunfield), or by Button~\cite{button2005}.
\end{rmk}

\begin{table}[htb]
\caption{Taut angle structures corresponding to veering triangulations for manifolds that do not fibre. The edges with angle $\pi$ are listed. The edge subscripts are the vertex numbers as given in the SnapPea census.\newline}
\centering
\begin{tabular}{ |r|ccccccc| }   
\hline
Tetrahedron  & 0 & 1 &2&3&4&5&6 \\ 
\hline
s227 & $e_{02}, e_{13}$ & $e_{01}, e_{23}$ & $e_{01}, e_{23}$ & $e_{02}, e_{13}$ & $e_{03}, e_{12}$ & $e_{03}, e_{12}$ & - \\
s438 & $e_{02}, e_{13}$ & $e_{02}, e_{13}$ & $e_{02}, e_{13}$ & $e_{02}, e_{13}$ & $e_{02}, e_{13}$ & $e_{01}, e_{23}$ & - \\ 
s772 & $e_{02}, e_{13}$ & $e_{01}, e_{23}$ & $e_{01}, e_{23}$ & $e_{01}, e_{23}$ & $e_{01}, e_{23}$ & $e_{02}, e_{13}$ & - \\
s773 & $e_{02}, e_{13}$ & $e_{01}, e_{23}$ & $e_{01}, e_{23}$ & $e_{01}, e_{23}$ & $e_{01}, e_{23}$ & $e_{02}, e_{13}$ & - \\
s779 & $e_{02}, e_{13}$ & $e_{01}, e_{23}$ & $e_{01}, e_{23}$ & $e_{01}, e_{23}$ & $e_{01}, e_{23}$ & $e_{02}, e_{13}$ & - \\
v3128 & $e_{02}, e_{13}$ & $e_{02}, e_{13}$ & $e_{02}, e_{13}$ & $e_{02}, e_{13}$ & $e_{02}, e_{13}$ & $e_{02}, e_{13}$ & $e_{01}, e_{23}$\\
v3243 & $e_{01}, e_{23}$ & $e_{02}, e_{13}$ & $e_{02}, e_{13}$ & $e_{02}, e_{13}$ & $e_{02}, e_{13}$ & $e_{01}, e_{23}$ & $e_{02}, e_{13}$ \\
v3244 & $e_{01}, e_{23}$ & $e_{03}, e_{12}$ & $e_{03}, e_{12}$ & $e_{03}, e_{12}$ & $e_{03}, e_{12}$ & $e_{01}, e_{23}$ & $e_{03}, e_{12}$ \\
v3377 & $e_{02}, e_{13}$ & $e_{02}, e_{13}$ & $e_{02}, e_{13}$ & $e_{02}, e_{13}$ & $e_{02}, e_{13}$ & $e_{01}, e_{23}$ & $e_{02}, e_{13}$ \\
v3526 & $e_{01}, e_{23}$ & $e_{03}, e_{12}$ & $e_{03}, e_{12}$ & $e_{02}, e_{13}$ & $e_{02}, e_{13}$ & $e_{03}, e_{12}$ & $e_{01}, e_{23}$ \\
\hline
\end{tabular}
\label{table_example_taut_angles}
\end{table}

\begin{rmk}\label{rmk:veering data}
There are 4,815 orientable triangulations in the SnapPea census, which lists triangulations with up to 7 tetrahedra. On those triangulations we calculate a total of 13,599 taut angle structures, of which 10,204 have a compatible taut structure and 158 have a compatible veering structure. There are 125 taut angle structures that are both taut and veering; all but one (s227) of the associated manifolds have some triangulation that is layered, according to Dunfield's data. Thus, we have at least $34$ veering structures not coming from Agol's construction.
\end{rmk}


\bibliography{}

\end{document}